\documentclass[12pt,oneside]{article}
\usepackage{amsmath,amssymb,amsfonts,amsthm}
\usepackage{amscd}
\usepackage{graphicx}
\usepackage[all]{xy}
\usepackage{multirow}
\usepackage{rotating}
\usepackage{hyperref}
\usepackage{hyperref}
\hypersetup{
    colorlinks=true,
    linkcolor=blue,
    filecolor=magenta,
    urlcolor=cyan,
    citecolor=blue}
\usepackage{lscape}
\newtheorem{Thm}{Theorem}[section]
\newtheorem{Lem}[Thm]{Lemma}

\newtheorem{Cor}[Thm]{Corollary}

\newtheorem{Def}[Thm]{Definition}
\newtheorem{Rem}[Thm]{Remark}

\theoremstyle{definition}

\theoremstyle{remark}

\renewcommand{\phi}{\varphi}

\newcommand{\Vol}{\operatorname{Vol}}

\renewcommand{\cos}{\operatorname{cos}}

\numberwithin{equation}{section}
\begin{document}
\title{Harmonic Finsler manifolds of $(\alpha , \beta)$-type}
%\vspace{-1.5 cm}
\author{Ebtsam  H. Taha}
\date{}
\maketitle    
\begin{center}
Department of Mathematics, Faculty of Science, Cairo University, Giza, Egypt,
\end{center}
\vspace{-0.8cm}
\begin{center}
E-mail: ebtsam.taha@sci.cu.edu.eg
\end{center}

\maketitle 

\vspace{0.5cm}
\textbf{Abstract.} In this paper we construct a new class of  harmonic  and asymptotically harmonic Finsler manifolds of $(\alpha , \beta)$-type.  
  This class is defined by a  Riemannian metric $\alpha$ and a special $1$-form $\beta$.\\
%\end{abstract}

%\tableofcontents
\maketitle 

\textbf{Keywords:}
Finsler metric; $(\alpha , \beta)$-metric; harmonic manifold; asymptotically harmonic;   Finsler mean curvature; Busemann-Hausdorff volume measure;  Holmes-Thompson volume measure.\\

\maketitle 
\textbf{MSC 2020:} 53B40, 53C60, 58B20, 58J60.

\section{Introduction}
\par A complete Riemannian manifold $(M, \alpha)$ is said to be harmonic if  for any $x \in M$, its volume density function $\sqrt{\det(\alpha_{ij}(x))}$ is a radial function (i.e., it depends only on the distance to the origin in the polar coordinates centered at $x$). All Riemannian harmonic manifolds $(M, \alpha)$ are Einstein manifolds, that is they have constant  Ricci curvature \cite{Besse}.  Further,  both  $2$-dimensional and $3$-dimensional harmonic Riemannian manifolds have constant sectional curvature.  Known examples of harmonic Riemannian manifolds are  the Euclidean space,  rank
one symmetric spaces  and  Damek–Ricci spaces \cite[Chapter 5]{harmonicbook16}, \cite{Besse, harmonicRes}.  

\vspace{5 pt}
%Studying harmonic geometry  is an active research area and 
It is known that, the density function of a Riemannian harmonic manifold plays a vital role  in studying many problems such as the classification of such manifolds \cite{Besse, harmonicRes}.  Finsler metrics are  natural generalization of Riemannian metrics. In \cite{harmonicFinsler}, harmonic manifolds have been extended to Finsler geometry. There are significant differences between the two geometries which are reflected in the study of harmonic Finsler spaces, such as the  dependence of the metric on  fiber coordinates, nonsymmetry of the Finsler distance, nonexistence of canonical  volume form, nonlinearity of the Finsler Laplacian, and appearance of non-Riemannian quantities. 

\vspace{5 pt}
Matsumoto in \cite{Matsumoto} gave an elegant method to build Finsler metrics from the Riemannian ones using a $1$-form. These metrics are known as $(\alpha , \beta)$-metrics and include Randers, Kropina, Matsumoto and square metrics for which many interesting results have been obtained cf. \cite{dV, harmonicFinsler, Shenbook16, arFinsler}.
 
 \vspace{5 pt}
In the present paper, we use the known machinery  of harmonic Riemannian manifolds to construct a new class of harmonic Finsler manifolds $(M, F, d \mu)$, namely \emph{harmonic Finsler manifolds of $(\alpha , \beta)$-type}, using either  Busemann-Hausdorff  volume measure or Holmes-Thompson volume measure on $M$. The main results of the present work are as follows: \\
\textbf{(a)} Theorem \ref{harmonic Kropina with constant length of beta} in which we construct a harmonic Kropina Finsler  manifold  $(M,F=\frac{\alpha^2}{\beta},\, d\mu)$ from a harmonic Riemmanian manifold  $(M,\alpha)$ and  a constant Killing $1$-form $\beta$ with respect to $\alpha$. Also, we prove that $(M,F:=\frac{\alpha^2}{\beta},\, d\mu)$ is of Einstein type. \\
\textbf{(b)} Theorem \ref{a thm} in which we provide, from a harmonic Riemmanian manifold $(M,\alpha)$ and  a $1$-form $\beta$ whose length $||\beta||_{\alpha}$ is a radial function, a harmonic  $(\alpha, \beta)$-Finsler manifold $(M,F=\alpha\, \phi(\frac{\beta}{\alpha}),\, d\mu)$.\\%, where $d\mu$ is either  Busemann-Hausdorff measure or Holmes-Thompson  volume measure.\\
 \textbf{(c)} It is known that, if $(M,\alpha)$ is noncompact harmonic Riemannian then it is asymptotically harmonic,  that is, $M$ has no conjugate points and the mean curvature of its horospheres is constant. We make use of this  observation and our above mentioned results to provide examples of asymptotically harmonic Finsler manifolds of  $(\alpha , \beta)$-type in Theorem \ref{AH thm}.
%This paper is organized as follows.  In Section 2, we recall some necessary materials for better understanding of our work. Section 3 is devoted to the main results in which we establish harmonic Finsler spaces of $(\alpha , \beta)$-type. Finally, we conclude some results  about asymptotically harmonic Finsler spaces of $(\alpha , \beta)$-type. 
\section{A Finsler geometric setting}
Let $M$ be a smooth connected $n$-dimensional manifold of dimension $n \geq 2$.  
 
 Let $(M,F)$ be a Finsler manifold. The Finsler distance $d_{F}$ induced by the Finsler structure $F$ \cite{Tam08} is defined  on $M$ by $$d_{F}(p,q):= \inf \left\{\int^{1}_{0} F(\dot{\xi}(t))\, dt\,\ | \, \xi : [0,1]  \rightarrow M,\,  C^{1}\text{ curve joining }p \text{ to } q \right\}.$$

 It should be noted that the \emph{Finsler distance} is nonsymmetric, that is, $d_{F}(p,q)\neq d_{F}(q,p)$.  In other words, the distance depends on the direction of the curve.  The non-reversibility property is also reflected in the notions of Cauchy sequence and completeness. Thus, being different from the Riemannian case,  a positively (or forward) complete Finsler manifold $(M , F)$ is not necessarily negatively (or backward) complete. For example, a Randers metric is forward complete solely.  A Finsler metric is said to be \emph{complete} if it is both forward and backward complete.

\vspace{5 pt}
It is worth mentioning that every function $\eta: \mathbb{R}^{+} \longrightarrow \mathbb{R}$ generates a \textit{radial function} $\eta_{x}$  around a point $x \in M$ defined by $\eta_{x}(x_{o})= \eta(r_{x}(x,x_{o}))$, where $r_{x}(x,x_{o})$ is the geodesic distance between $x, x_{o} \in M$ induced by the Finsler function \cite{harmonicFinsler}.
 %M. Matsumoto, in \cite{Matsumoto}, introduced the celebrated $(\alpha, \beta)$-metrics. 
 \begin{Def}\emph{\cite{Matsumoto}}
Let $\alpha = \sqrt{\alpha_{ij}(x)\, y^{i}\,y^{j}}$ be a Riemannian metric defined on $M$ and $\beta = b_{i}(x)\, dx^{i}$  be a $1$-form on $M$ with $b:=||\beta||_{\alpha} < b_{o}$,  where $b_{o}$ is a real positive number.  The Finsler metric $F=\alpha\, \phi(s),\, s:=\frac{\beta}{\alpha},$  is said to be an $(\alpha, \beta)$-metric if $\phi$ is a positive smooth  function defined on the interval $(-b_{o}, b_{o})$ such that \[\phi(s)- s\, \phi'(s) +(b^2 -s^2 ) \,\phi''(s) > 0,\,\,\,\, |s| \leq b < b_{o}. \]
\end{Def}
%Its metric tensor is given by {\small{ \begin{equation}\label{gen metric} g_{ij}= \left(\phi^2 -s \,\phi\, \phi' \right)\, \alpha_{ij} + \left( \phi \, \phi'' + (\phi')^{2}\right)\, b_{i}\, b_{j}+ \frac{1}{\alpha}\left(\phi\, \phi' -s\left[\phi\,\phi'' + (\phi')^{2} \right] \right) \, \left(b_{i} y_{j} + b_{j} y_{i} -  \frac{s}{\alpha} y_{i} y_{j}\right),\end{equation} }}where $y_{j}:= \alpha_{ij}\, y^{i}$.\\
Unlike the Riemannian case, there are several  volume forms used in Finsler geometry.  The most important ones are \emph{Busemann-Hausdorff and Holmes-Thompson volume forms} \cite{Shenbook16}.

%Let $(M,F)$ be n-dimensional smooth orientable manifold. 
\begin{Def}
The Busemann-Hausdorff volume form $d\mu_{BH}:= \sigma_{BH}(x)\,dx$ is defined  at a point $x \in M$ as follows 
\begin{equation*}
\sigma_{BH}(x):= \frac{\Vol(\mathbb{B}^n(1))}{\Vol(B_{F}^n(1))},\,\,\quad\text{where}\,\, \Vol(B_{F}^n(1)):= \Vol \left(\{ (y^i) \in \mathbb{R}^{n}\,| \,F(x,y^i \partial_{i})<1\}\right)
\end{equation*}
and $\Vol (\mathbb{B}^n(1))$ denotes the Euclidean volume of a unit Euclidean ball:
$$\Vol(\mathbb{B}^n(1))=\frac{1}{n} \Vol(\mathbb{S}^{n-1})= \frac{1}{n} \Vol(\mathbb{S}^{n-2}) \int_{0}^{\pi} \sin^{n-2}t\,dt. $$
\end{Def}
%For a general Finsler manifold, the Busemann-Hausdorff volume form may be expressed hardly by element functions (the Finsler structure). However, it was done for particular Finsler metrics, namely Randers metrics \cite{Shenbook16, volform}.
\begin{Def}The {Holmes-Thompson volume form}  is defined  at a point $x \in M$ by $d\mu_{HT}: = \sigma_{HT}(x) \,dx,$ where
$$\sigma_{HT}(x):= \frac{1}{\Vol(\mathbb{S}^{n-1})} \int_{S_{x}M} \det(g_{ij}(x,y))\, dy,$$ where 
$S_{x}M = \{ y \in T_{x}M\,| \,F(y)=1\}$.
%$$\sigma_{HT}(x):=\frac{\Vol(B_{F}^n(1), g)} {\Vol(\mathbb{B}^n(1))}= \frac{1}{\Vol(\mathbb{B}^n(1))} \int_{B_{F}^n(1)} \det(g_{ij}(x,y))\, dy.$$
\end{Def}

\begin{Def}\emph{\cite{harmonicFinsler}}\label{harF}
\it{A forward complete Finsler manifold $(M, F,  d\mu)$ endowed  with a smooth volume measure $d\mu$  is (globally) harmonic  if in polar coordinates the volume density function $\overline{\sigma}_{p}(r,y)$ is a radial function around (each) $ p \in M,$ where $\overline{\sigma}_{p}(r,y):= \frac{{\sigma}_{p}(r,y)}{\sqrt{\det(\dot{g}_p(p,y))}}$  and $\dot{g}_p$ is the restriction of $g$ on the indricatrix $I_{p}M$.  That is,  $\overline{\sigma}_{p}(r,y)$ is independent of $y \in I_{p}M$; so it can be written briefly as  $\overline{\sigma}_{p}(r)$.}
\end{Def}
There are some equivalent definitions to Definition \ref{harF} that can be found in \cite{harmonicFinsler}.
\begin{Def}\emph{\cite{harmonicFinsler}}
Suppose that $r$ is a Finsler distance defined on an open subset of  $M$, that is $F(\nabla r)=1$, where $\nabla r$ is the gradient of $r$. The Finsler mean curvature of the level hypersurface $r^{-1}(t)$ at $x \in M$  with respect to $\nabla r_{x}$ is defined by  
\begin{equation}\label{Fmeancurvaturedef}
\Pi_{\nabla r}(x):= \frac{d}{dt}\,\log(\sigma_{x}(t,x^{a}))|_{t=t_{o}},
\vspace{-0.5 cm}
\end{equation}
for some  $t_{o} \in$  \emph{Im}$(r)$.
\end{Def}
\begin{Def}\emph{\cite{harmonicFinsler}}
A forward complete, simply connected Finsler manifold $(M, F,d \mu)$ without conjugate points  is called an \emph{asymptotically harmonic Finsler manifold (or simply, AHF-manifold)} if the  mean curvature of Finslerian horospheres (Finsler geodesic sphere with infinite radius) $\Pi_{\infty}:=\lim_{r \rightarrow \infty}\Pi_{\nabla r}$ is a real constant.
\end{Def}
 Consequently, a simply connected, forward complete, noncompact harmonic Finsler manifold with  constant Finsler mean curvature of horospheres is an AHF-manifold cf. \cite[Chapter 5]{harmonicbook16}, \cite{ harmonicFinsler}.
 
For general Finsler metrics, the Busemann-Hausdorff  and Holmes-Thompson  volume forms may be expressed hardly by elementary functions. However, it is done for some particular Finsler metrics, namely $(\alpha , \beta)$-Finsler metrics. %More precisely,  it was proved in  that
\begin{Lem}\emph{\cite{dV}}\label{vol forms}
The Busemann-Hausdorff volume form of the $(\alpha , \beta)$-Finsler metric $F=\alpha\, \phi(s)$, with $||\beta||_{\alpha}=b$, is given by \vspace{-0.3 cm}
\[  d\mu_{BH}= f_{BH}(b)\,\, d\mu_{\alpha}, \qquad   f_{BH}(b)= \frac{ \int_{0}^{\pi} \sin^{n-2}t\, dt}{\int_{0}^{\pi} \frac{\sin^{n-2}t}{\phi^{n}(b\, \cos t)}\, dt} .\,\,\,\]
The Holmes-Thompson volume measure is given by
\[  d\mu_{HT}= f_{HT}(b)\,\, d\mu_{\alpha}, \qquad  f_{HT}(b)= \frac{\int_{0}^{\pi} T(b\, \cos t)\, \sin^{n-2}t \, dt}{\int_{0}^{\pi} \sin^{n-2}t\, dt},\,\,\,\] where $d\mu_{\alpha} =\sqrt{\det(\alpha_{ij}(x))}\, dx$ denotes the Riemannian volume form of $\alpha$ and $T$ is defined by 
\begin{equation}\label{T}
T(s)= \phi (\phi -s \phi')^{n-2} \left[(\phi - s\phi')+ (b^2 -s^2)\, \phi'' \right].
\end{equation}
\end{Lem}
\section{Harmonic and asymptotically harmonic Finsler manifolds of  $(\alpha , \beta)$-type}
%are the following:
%%%%%%%%%%%%%%%%%%%%%%%%%%%%%%%%%%%%%%%%%%%%
\begin{Def}\emph{\cite{Shenbook16}}
A $1$-form $\beta$ is called a constant Killing $1$-form with respect to $\alpha$ if $||\beta||_{\alpha}$ is constant and $r_{ij} :=\frac{1}{2}( b_{i|j} -b_{j|i})=0$, where $b_{i|j}$  denotes covariant derivative of $\beta$ with respect to $\alpha$.
\end{Def}
 %Apparently, in Riemannian manifold $(M, \alpha)$: the Einstein condition ($Ric_{\alpha} = \lambda\, \alpha(x)$ for some constant $\lambda$) implies that harmonic Riemannian manifolds belong to one of the following three classes cf. \cite[Ch.5]{harmonicbook16}: (1) positive Ricci curvature $Ric_{\alpha}  > 0$ implies $(M, \alpha)$ is a compact rank one symmetric space which results from Myers-Bonnet theorem, (2) vanishing Ricci curvature  $ Ric_{\alpha}  = 0$ implies $(M, \alpha)$ is the Euclidean space $\mathbb{R}^n$ which results by the Riccati equation and (3) negative Ricci curvature $ Ric_{\alpha} < 0$ implies $(M, \alpha)$ is noncompact and Damek-Ricci space.\\

Hereafter, we assume that $M$ is a smooth connected forward complete manifold and $d\mu$ is either  Busemann-Hausdorff or Holmes-Thompson  volume measures on $M$.

\begin{Thm}\label{harmonic Kropina with constant length of beta}
Let $(M,\alpha)$ be a harmonic Riemmanian manifold of dimension $n \geq 2$. Let $\beta$ be a constant Killing $1$-form with respect to $\alpha$ such that  $||\beta||_{\alpha} < b_{o}$. Then $(M,F=\frac{\alpha^2}{\beta},\, d\mu)$ is  harmonic non-Riemanian Finsler manifold of Einstein type.
\end{Thm}
\begin{proof}
As $(M,\alpha)$ is a harmonic Riemmanian manifold, the volume density function $\sqrt{\det(\alpha_{ij})}$ of $d\mu_{\alpha}$ is a radial function, say $l(r)$. In other words, $d\mu_{\alpha}= l(r)\, dr \wedge d\Theta$.
Since $\beta$ is a constant Killing $1$-form with respect to $\alpha$, thereby $||\beta||_{\alpha}$ is constant, say $A$.  Then, in view of Lemma \ref{vol forms}, we have  for $\phi(s)= \frac{1}{s}, \, s >0$,
\[  d\mu_{BH}= \frac{ \int_{0}^{\pi} \sin^{n-2}t\, dt}{\int_{0}^{\pi} \frac{\sin^{n-2}t}{\phi^{n}(A\, \cos t)}\, dt} \,\,l(r)\, dr \wedge d\Theta\,\]
and 
\[  d\mu_{HT}=  \frac{\int_{0}^{\pi} T(A\, \cos t)\,\sin^{n-2}t \, dt}{\int_{0}^{\pi} \sin^{n-2}t\, dt}\,\,l(r)\, dr \wedge d\Theta.\]
Hence, the corresponding volume density functions are radial functions. The Kropina Finsler manifold $(M,F=\frac{\alpha^2}{\beta}, d\mu)$ is a harmonic since both  $d\mu_{BH}$ and $d\mu_{HT}$ are constant multiples of $d\mu_{\alpha}$, which is radial. Moreover, it is of Einstein type in view of \cite[Theorem 10.8]{Shenbook16} $\,(\alpha$ is Einstein; being harmonic Riemannian, and $\beta$ is constant killing).
\end{proof}
\begin{Rem}\emph{
The Einstein harmonic Kropina Finsler manifold $(M,F=\frac{\alpha^2}{\beta},\, d\mu)$ appearing in Theorem \ref{harmonic Kropina with constant length of beta} has the properties that:\\
$(1)$ it has constant Einstein factor of $F$. This  follows from \cite[Theorem 10.8]{Shenbook16}.\\
$(2)$ it  has constant flag curvature for $n=2,\, n=3$. This follows from the fact that all harmonic Riemannian spaces of dimension 2 and 3 have constant sectional curvature  \cite{Besse} together with \cite[Corollary 10.3]{Shenbook16}.}
\end{Rem}
%t is known that, if $M$ is noncompact and harmonic then it is asymptotically harmonic,  that is, $M$ has no conjugate points and the mean curvature of its horospheres is constant. 
%\begin{Thm}\label{harmonic(alpha,beta)thm}Let $(M,\alpha)$ be a harmonic Riemmanian manifold. Let $\beta$ be a $1$-form such that its length $||\beta||_{\alpha}$ is a radial function and $||\beta||_{\alpha} < 1$. Then $(M,F:=\alpha\, \phi \left(\frac{\beta}{\alpha}\right),\, \mu)$ is harmonic Finsler manifold, where $d\mu$ is either  Busemann-Hausdorff or Holmes-Thompson  volume measures on $M$.\end{Thm}
Now, we provide a constructive way to find a new class of harmonic  Finsler manifolds.
\begin{Thm}\label{dVHT harmonic(alpha,beta)thm}
Let $(M,\alpha)$ be a harmonic Riemmanian manifold. Let $ T(s)$ be defined by \eqref{T}. If  
$\,T(s)-1$ is an odd function of $s$, then $(M,F=\alpha\, \phi (s),\, d\mu_{HT})$ is a harmonic  Finsler manifold. % where $d\mu$ is  Holmes-Thompson  volume measure on $M$.
\end{Thm}
\begin{proof}
As $(M,\alpha)$ is a harmonic Riemmanian manifold, the volume density function $\sqrt{\det(\alpha_{ij})}$ of $d\mu_{\alpha}$ is a radial function, say $l(r)$, that is, $d\mu_{\alpha}= l(r)\, dr \wedge d\Theta$. 
 Now,  if
$T(s)-1$ is an odd function of $s$, then, according to  \cite[Corollary 2.2]{dV}, the Holmes-Thompson volume measure coincides with  $d\mu_{\alpha}$. That is,
$$ d\mu_{HT}=  d\mu_{\alpha} = l(r)\, dr \wedge d\Theta.$$
Hence, $(M,F=\alpha\, \phi (s),\, d\mu_{HT})$ is a harmonic  Finsler manifold of $(\alpha , \beta)$-type.
\end{proof}

It is worth mentioning that $T(s)-1$ is an odd function of $s$ for some significant particular $\phi$ and thereby Theorem \ref{dVHT harmonic(alpha,beta)thm} can be applied. For example, in the case of Randers metrics, $\phi (s) = 1+s$ and by  \eqref{T}, $T(s)= 1+s$, which implies that $T(s)-1= s$ is an odd function of $s$. Hence, we have the result%Further, the harmonic Finsler manifolds of Randers type  with respect to both $d\mu_{BH},\,d\mu_{HT}$ and extreme volume forms have been studied in \cite[\S 5]{harmonicFinsler} using a different technique.
\begin{Cor}
All Randers manifolds $(M, F= \alpha + \beta, d\mu_{HT})$ with harmonic $\alpha$ are harmonic.
\end{Cor}
%$(M, F= \alpha + \beta, d\mu_{HT})$
\begin{Thm}\label{a thm}
Let $(M,\alpha)$ be a harmonic Riemmanian manifold. Let $\beta$ be a $1$-form such that its length $||\beta||_{\alpha}$ is a radial function and $||\beta||_{\alpha} < b_o$. Then, $(M,F=\alpha\, \phi(s),\, d\mu)$ is a harmonic Finsler manifold of $(\alpha, \beta)$-type.
\end{Thm}
\begin{proof}
As $(M,\alpha)$ is a harmonic Riemmanian manifold, the volume density function $\sqrt{\det(\alpha_{ij})}$ of $(M,{\alpha})$ is a radial function, say $l(r)$. In other words, $d\mu_{\alpha}= l(r)\, dr \wedge d\Theta$.
Since,  $||\beta||_{\alpha}$  is a radial function.  Then, in view of Lemma \ref{vol forms}, we obtain 
\[  \sigma_{BH}= f_{BH}(r)\,l(r) \qquad \text{ and }\, \qquad \sigma_{HT}=  f_{HT}(r)\,l(r).\]
%and the Holmes-Thompson volume measure is given by\[  \sigma_{HT}=  f_{HT}(r)\,l(r)\, dr \wedge d\Theta.\]
Hence, the corresponding volume density functions are radial functions. Consequently, $(M,F=\alpha\phi(s), d\mu)$ is  harmonic  Finsler manifold of $(\alpha, \beta)$-type.
\end{proof}

We end this work by applying our aforementioned results to asymptotically harmonic Finsler manifolds of  $(\alpha , \beta)$-type. The mean curvature of a Finslerain geodesic sphere is denoted by  $\Pi^{F}_{\nabla r}$ and the mean curvature of a Riemannian  geodesic sphere is denoted by  $\Pi^{\alpha}_{\nabla r}$. Consequently,  the mean curvature of Finslerian horospheres and Riemannian horospheres are denoted by  $\Pi^{F}_{\infty},\, \Pi^{\alpha}_{\infty}$, respectively.

\begin{Thm}\label{AH thm}
Let $(M,\alpha)$ be a complete, simply connected noncompact harmonic Riemmanian space with horospheres of constant mean curvature. The following assertions hold:\\
$(1)$ If $\beta$ is a constant Killing $1$-form with respect to $\alpha$, then $(M,F=\frac{\alpha^2}{\beta},\, d\mu)$ is an AHF-manifold.\\
$(2)$ If $T(s)-1$ is an odd function of $s$, where $T(s)$ is defined by \eqref{T}, then $(M,F=\alpha\, \phi(s),\, d\mu_{HT})$ is an AHF-manifold.\\
$(3)$ If  $||\beta||_{\alpha}$ is radial,  then $(M,F=\alpha\, \phi(s),\, d\mu)$ is an AHF-manifold.
\end{Thm}
\begin{proof}
Suppose $(M,\alpha)$ is a  harmonic Riemmanian manifold. \\
$(1)$ By Theorem \ref{harmonic Kropina with constant length of beta}, $(M,F=\alpha^2 /\beta,\, d\mu)$ is a harmonic Finsler manifold. Moreover, each of $ d\mu_{HT}$ and $d\mu_{BH}$ is a constant multiple of $d\mu_{\alpha}$. Hence, the mean curvature of geodesic sphere $\Pi^{F}_{\nabla r}$ is equal to  $\Pi^{\alpha}_{\nabla r}$, by \eqref{Fmeancurvaturedef}. Now, since $(M,\alpha)$ is a complete, simply connected noncompact harmonic Riemmanian manifold with horospheres of constant mean curvature, that is,  $\Pi^{\alpha}_{\infty}$ exists.  Consequently,  $\Pi^{F}_{\infty}$ is a real constant. Then it is an AHF-manifold. \\
$(2)$ Assume $T(s)-1$ is an odd function of $s$. By Theorem \ref{dVHT harmonic(alpha,beta)thm}, $(M,F=\alpha\, \phi(s),\, d\mu_{HT})$ is a harmonic  Finsler manifold and $ d\mu_{HT} = d\mu_{\alpha}$. Thus, by \eqref{Fmeancurvaturedef}, $\Pi^{F}_{\nabla r} = \Pi^{\alpha}_{\nabla r}$. Since $(M,\alpha)$ is a complete, simply connected noncompact harmonic Riemmanian manifold with horospheres of constant mean curvature, that is,  $\Pi^{\alpha}_{\infty}$ exists; thereby,  $\Pi^{F}_{\infty}$ exists. Hence, $(M,F=\alpha\, \phi (s),\, d\mu_{HT})$ is an AHF-manifold.\\ \textbf{$(3)$} It follows directly from Theorem \ref{a thm}.
\end{proof}
%\vspace{-1.3 cm}

\textbf{Acknowledgment.} I would like to express my deep thanks to Professor Nabil L. Youssef for his careful reading of this manuscript and his valuable discussions which led to the present version.

%%%%%%%%%%%%%%%%%%%%%%%%%%%%%%%%%%%%%%%


\begin{thebibliography}{}

\bibitem{harmonicbook16}C. S. Aravinda, F. T. Farrell and  J.-F. Lafont (Eds.), \textit{Geometry, Topology, and Dynamics in Negative Curvature} (London Mathematical Society Lecture Note Series), Cambridge University Press, 2016. 

\bibitem{Besse}A. L. Besse,
\emph{Manifolds all of whose geodesics are closed}, Springer-Verlag, Berlin-New York, 1978. %ix+262 pp. ISBN: 3-540-08158-5.

\bibitem{dV} X. Cheng and Z. Shen,\textit{ A class of Finsler metrics with isotropic S-curvature},  Israel J.  Math.,  \textbf{169} (2009), 317-340.


\bibitem{harmonicRes} P. B. Gilkey and J. H. Park, \textit{Harmonic spaces and density functions}, Results Math., \textbf{75}:121 (2020). 


%\bibitem{Shlec}Z.~Shen, \emph{Lectures on Finsler geometry},World Scientific Publishing Co., Singapore, 2001.

\bibitem{Matsumoto} M. Matsumoto, \textit{On C-reducible Finsler spaces}, Tensor N. S., \textbf{24} (1972), 29-37.

\bibitem{harmonicFinsler} H.  Shah and E.  H. Taha, \emph{On Harmonic and asymptotically harmonic Finsler manifolds}, Submitted, arXiv: 2005.03616v3 [math.DG].



\bibitem{Shenbook16} Y. Shen and Z. Shen, \emph{Introduction to modern Finsler geometry}, World Scientific Publishing Co., Singapore, 2016.% xii+393 pp. ISBN: 978-981-4704-90-8.

 
 \bibitem{arFinsler} E. H. Taha and B. Tiwari, \emph{On almost rational Finsler metrics}, Submitted, arXiv: 2101.01764 [math.DG].
 
 \bibitem{Tam08}  L. Tam\'assy, \emph{Distance functions of Finsler spaces and distance spaces,} Proceeding of the conference \lq \lq Differential  Geometry and its  Applications, 2008", pp. 559-570. World Scientific Publ., 2008.


%\bibitem{Szabo} Z. I. Szab\'o, \emph{The Lichnerowicz conjecture on harmonic manifolds} J. Differential Geom. 31 (1990), no. 1, 1-28.

\end{thebibliography}
\end{document}